\documentclass[12pt]{amsart}

\usepackage{amsmath,amssymb,amsthm}
\usepackage{epsf}
\usepackage{epsfig}
\usepackage[T2A]{fontenc}
\usepackage[cp1251]{inputenc}
\usepackage[english]{babel}

\newcommand{\RR}{\mathbb R}
\newcommand{\CC}{\mathbb C}

\renewcommand{\Re}{\mathop{\rm Re}\nolimits}
\renewcommand{\Im}{\mathop{\rm Im}\nolimits}

\newcommand{\ii}{\mathrm{i}}

\newcommand{\manifold}[1]{\mathcal{#1}}
\newcommand{\M}{\manifold{M}}

\newcommand{\vect}[1]{\mathrm{#1}} 
\newcommand{\x}{\vect{x}}

\newcommand{\vb}{\vect{b}}

\newcommand{\vX}{\vect{X}}
\newcommand{\vY}{\vect{Y}}

\newtheorem{thm}{Theorem}[section]

\theoremstyle{definition}

\theoremstyle{remark}

\newcommand{\prfend}{\hfill $\blacksquare$ \bigskip}

\newcommand{\ds}{\displaystyle}

\makeatletter \@addtoreset{equation}{section} \makeatother

\begin{document}

\title[EXPLICIT SOLVING OF THE SYSTEM OF NATURAL PDE'S]
{Explicit Solving of the System of Natural PDE's\\ of Minimal Space-like Surfaces\\ in Minkowski Space-time}%

\author{Georgi Ganchev and Krasimir Kanchev}

\address{Bulgarian Academy of Sciences, Institute of Mathematics and Informatics,
Acad. G. Bonchev Str. bl. 8, 1113 Sofia, Bulgaria}
\email{ganchev@math.bas.bg}%

\address {Department of Mathematics and Informatics, Todor Kableshkov University of Transport,
158 Geo Milev Str., 1574 Sofia, Bulgaria}%
\email{kbkanchev@yahoo.com}%


\begin{abstract}

A minimal space-like surface in Minkowski space-time is said to be of general type
if it is free of degenerate points. The fact that minimal space-like surfaces of general type
in Minkowski space-time admit canonical parameters of the first (second) type implies that
any minimal space-like surface is determined uniquely up to a motion in space-time by the Gauss
curvature and the normal curvature, satisfying a system of two PDE's (the system of natural PDE's).
In fact this solves the problem of Lund-Regge for minimal space-like surfaces. Using
canonical Weierstrass representations of minimal space-like surfaces of general type in
Minkowski space-time we solve explicitly the system of natural PDE's, expressing any solution
by means of two holomorphic functions in the Gauss plane. We find the relation between two pairs of
holomorphic functions (i.e. the class of pairs of holomorphic functions) generating
one and the same solution to the system of natural PDE's, i.e. generating one and the same minimal
space-like surface in Minkowski space-time.

\end{abstract}

\maketitle

\thispagestyle{empty}

\section{Introduction}
We study space-like surfaces $\M$ in Minkowski space-time $\RR^4_1$, i.e. surfaces whose induced metric
is of signature $(2,0)$. The surface $\M$ is minimal if its mean curvature vector field $H$ is zero.

Let $(\M,\x(u,v))$ be a space-like surface in $\RR^4_1$, parameterized by isothermal coordinates $(u, v)$.
We denote by $T_p(\M)$ and $N_p(\M)$, the tangential space and normal space
at a point $p \in \M$, respectively. The second  fundamental form on $\M$ is denoted by $\sigma$.

\emph{A point $p\in\M$ is said to be degenerate, if the set $\{\sigma(\vX,\vY);\ \vX \in T_p(\M),\vY \in T_p(\M) \}$,
is contained in one of the two light-like one-dimensional subspaces of $N_p(\M)$.}

A point $p \in \M$ is degenerate if and only if the Gauss curvature $K$ and the curvature of the normal connection
$\varkappa$ (the normal curvature) are zero at the point $p$.

\emph{We call a minimal space-like surface, free of degenerate points, a
minimal space-like surface of general type.}

The Gauss curvature $K$ and the normal curvature $\varkappa$
of a minimal space-like surface of general type, parameterized by special isothermal parameters, satisfy
the following system of  partial differential equations \cite{A-P-1}:
\begin{equation}\label{Mink-SL-NatEq}
\begin{array}{l}
\ds{(K^2 + \varkappa^2)^{\frac{1}{4}}\, \Delta \ln (K^2+\varkappa^2)^\frac{1}{4}} = 2K,\\
[2mm]
\ds{(K^2 + \varkappa^2)^{\frac{1}{4}}\, \Delta \arctan
\frac{\varkappa}{K} = 2\varkappa}.
\end{array}
\end{equation}

Conversely, any solution ($K$, $\varkappa$) to system \eqref{Mink-SL-NatEq}
determines uniquely (up to a motion in $\RR^4_1$) a minimal space-like surface of general type
with Gauss curvature $K$ and normal curvature $\varkappa$.

Further we call system \eqref{Mink-SL-NatEq} {\it the system of natural PDE's of minimal space-like surfaces
in $\RR^4_1$} and our aim is to solve explicitly this system.

All considerations in the paper are local.

In \cite{G-K-3} we proved that any minimal space-like surface of general type admits locally
canonical parameters of the first (second) type. The special isothermal parameters in \cite{A-P-1} appear
to be canonical parameters of the first type. These parameters are characterized by the second
fundamental form in the following way:
\begin{equation}\label{can1}
\begin{array}{l}
\sigma (\x_u,\x_u)\bot \: \sigma (\x_u,\x_v),\\[2mm]
\sigma^2(\x_u,\x_u)-\sigma^2(\x_u,\x_v)=1.
\end{array}
\end{equation}

If we introduce canonical parameters on any minimal space-like surface of general type,
then the Gauss curvature $K$ and the normal curvature $\varkappa$ satisfy the system of natural equations
\eqref{Mink-SL-NatEq} and determine the surface up to a motion in $\RR^4_1$. It is clear that the number
of the the invariants and the number of the PDE's can not be reduced further. Therefore this solves
the problem of Lund-Regge \cite{L-R} for minimal space-like surfaces in $\RR^4_1$.

In this paper we prove the following theorems.
\vskip 2mm
\textbf{Theorem 1.} {\it If the pair $(K,\varkappa)$ is a solution to the system \eqref{Mink-SL-NatEq}, then
the curvatures $K$ and $\varkappa$ are given locally by the formulas
\begin{equation}\label{Mink-SL-NatEq-Sol-g12}
\begin{array}{lll}
        K         &=& |\alpha|\Re\alpha \\
        \varkappa &=& |\alpha|\Im\alpha
\end{array}; \quad
\alpha = \ds\frac{-4g'_1 \bar {g'_2}}{(1 + g_1 \bar g_2)^2}\;,
\end{equation}
where $(g_1,g_2)$ are holomorphic functions satisfying the conditions
\begin{equation}\label{g_1g_2-cond}
g'_1 g'_2 \neq 0; \quad g_1 \bar g_2 \neq -1.
\end{equation}
Conversely, any pair of holomorphic functions satisfying the conditions \eqref{g_1g_2-cond} generates by means of
equalities \eqref{Mink-SL-NatEq-Sol-g12} a solution to the system \eqref{Mink-SL-NatEq}.}
\vskip 2mm
\textbf{Theorem 2.}{\it Two pairs of holomorphic functions $(g_1,g_2)$ and $(\hat g_1,\hat g_2)$
generate by means of \eqref{Mink-SL-NatEq-Sol-g12} one and the same solution to the system \eqref{Mink-SL-NatEq},
if and only if they are related by equalities of the type:
\begin{equation}\label{Mink-SL-Sol_hatg12_g12}
\hat g_1 = \ds\frac{ag_1+b}{cg_1+d}\,; \quad \hat g_2 = \ds\frac{\phantom{-}\bar d g_2 - \bar c}{-\bar b g_2 + \bar a}\;,
\end{equation}
where $a,b,c,d\in \CC$, $ad-bc\neq 0$.}
\vskip 2mm
\section{Explicit solving of the system of natural PDE's of minimal space-like surfaces of general type in $R_1^4$}

In \cite{G-M-1} it is shown that using appropriate substitutions the above system \eqref{Mink-SL-NatEq}
can be simplified. More precisely, putting $K=e^{2X}\cos Y$; $\varkappa=e^{2X}\sin Y$, we get the following
form of the natural equations:
\begin{equation}\label{Mink-SL-NatEq_XY}
    \begin{array}{lrl}
        \Delta X &=& 2e^{X}\cos Y,\\
        \Delta Y &=& 2e^{X}\sin Y.
    \end{array}
\end{equation}

It is easy to see that the system \eqref{Mink-SL-NatEq_XY} of natural equations is equivalent to one
complex partial differential equation:
\begin{equation}\label{Mink-SL-NatEq-Compl-1}
\Delta (X+iY) = 2e^{X}(\cos Y + i\sin Y) = 2e^{X+iY}.
\end{equation}
Putting in the last equation
\begin{equation}\label{alpha-SL-1}
\alpha = e^{X+iY},
\end{equation}
we have:
\begin{equation}\label{Mink-SL-NatEq-Compl-2}
\Delta \log\alpha = 2\alpha.
\end{equation}
The last equation coinsides formally with the natural partial differential equation (1.1) in \cite{G-2} of
space-like surfaces in $\RR_1^3$ with zero mean curvature:
\begin{equation}\label{Mink-SL-NatEq-R31}
\Delta \ln\nu - 2\nu = 0.
\end{equation}
The difference between \eqref{Mink-SL-NatEq-Compl-2} and \eqref{Mink-SL-NatEq-R31} is that $\nu$ is a real
positive function, while $\alpha$ is a complex function. According to Theorem 6.2 in \cite{G-2} any solution to
\eqref{Mink-SL-NatEq-R31} is given by the formula:
\begin{equation}\label{Mink-SL-NatEq-R31-Sol-1}
\nu (u,v) = \ds\frac{{\eta'_u}^2(u,v)+{\eta'_v}^2(u,v)}{\eta^2(u,v)},
\end{equation}
where $\eta(u,v)$ is an arbitrary real harmonic function, satisfying the condition: $\eta\neq 0$
and ${\eta'_u}^2+{\eta'_v}^2\neq 0$.
By a direct substitution of \eqref{Mink-SL-NatEq-R31-Sol-1} into
\eqref{Mink-SL-NatEq-Compl-2} it is easy to check that this formula gives a solution also in the case, when
$\eta(u,v)$ is an arbitrary complex harmonic function. Taking into account \eqref{alpha-SL-1} we use the following denotations:
\begin{equation}\label{eX-alpha}
e^{X} = |\alpha|; \quad e^{X}\cos Y = \Re\alpha; \quad e^{X}\sin Y
= \Im\alpha.
\end{equation}
Thus we found the following family of solutions to the system \eqref{Mink-SL-NatEq} of natural equations of minimal space-like
surfaces of general type in $\RR_1^4$:
\begin{equation}\label{Mink-SL-NatEq-Sol-1}
\begin{array}{lll}
        K         &=& |\alpha|\Re\alpha \\
        \varkappa &=& |\alpha|\Im\alpha
\end{array}; \quad
\alpha = \ds\frac{{\eta'_u}^2+{\eta'_v}^2}{\eta^2}\;,
\end{equation}
where $\eta(u,v)$ is an arbitrary complex harmonic function, satisfying the conditions: $\eta\neq 0$
and ${\eta'_u}^2+{\eta'_v}^2\neq 0$.

The following question arises naturally:

\emph{Does the formula \eqref{Mink-SL-NatEq-Sol-1} give all the solutions to the system
\eqref{Mink-SL-NatEq} of natural partial differential equations of minimal space-like surfaces in $\RR_1^4$?}

In \cite{G-K-3} we proved the following statement:
\begin{thm}
Any minimal space-like surface ${\M}$ of general type,
parameterized by canonical coordinates of the first type, has the
following Weierstrass representation:
\begin{equation}\label{Wcanh}
\Phi: \quad
\begin{array}{rlr}
\phi_1 &=& \ii\ds\frac{\cosh h_1}{\sqrt{{h'_1}^2 - {h'_2}^2}}\,,\\[8mm]
\phi_2 &=&  \ds\frac{\sinh h_1}{\sqrt{{h'_1}^2 - {h'_2}^2}}\,,\\[8mm]
\phi_3 &=&  \ds\frac{\cosh h_2}{\sqrt{{h'_1}^2 - {h'_2}^2}}\,,\\[8mm]
\phi_4 &=&  \ds\frac{\sinh h_2}{\sqrt{{h'_1}^2 - {h'_2}^2}}\,,\\
\end{array}
\end{equation}
where $(h_1,h_2)$ are holomorphic functions satisfying the
conditions:
\begin{equation}\label{Wcanh_cond}
{h'_1}^2 \neq {h'_2}^2; \quad \Re h_1\neq 0\ \text{or}\ \Im
h_2\neq \frac{\pi}{2}+k\pi;\ k\in {\mathbb Z}.
\end{equation}

Conversely, if $(h_1,h_2)$ is a pair of holomorphic functions
satisfying the conditions \eqref{Wcanh_cond}, then formulas
\eqref{Wcanh} generate a minimal space-like surface of general
type, parameterized by canonical coordinates of the first type.
\end{thm}
Introducing the functions:
$$w_1=h_1+h_2, \quad w_2=h_1-h_2; \quad g_1=e^{w_1}; \quad g_2=e^{w_2},$$
$$\theta=\Re h_1+\ii \Im h_2,$$
then \eqref{Wcanh} implies that \cite{G-K-3}:

\begin{equation}\label{K+ikappa_Can_theta}
K+\ii\varkappa  = \ds\frac{-|{\theta'_u}^2+{\theta'_v}^2|\, ({\theta'_u}^2+{\theta'_v}^2)}
                          {|\cosh \theta|^2\; \cosh^2 \theta}\;,
\end{equation}

\begin{equation}\label{K+ikappa_Can_w12}
K+\ii\varkappa  = \ds\frac{- |w'_1 w'_2|\, w'_1 \bar {w'_2}}
{\left|\cosh\frac{w_1 + \bar w_2}{2}\right|^2\; \cosh^2\frac{w_1 + \bar w_2}{2}}\;,
\end{equation}

\begin{equation}\label{K+ikappa_Can_g12}
K+\ii\varkappa  = \ds\frac{-16|g'_1 g'_2|\, g'_1 \bar {g'_2}}
                          {|1 + g_1 \bar g_2|^2\; (1 + g_1 \bar g_2)^2}\;.
\end{equation}

In terms of the functions $w_1, w_2$ the canonical Weierstrass representation of minimal
space-like surfaces of general type has the following form \cite{G-K-3}:
\begin{thm}
Any minimal space-like surface ${\M}$ of general type,
parameterized by canonical coordinates of the first type, has the
following Weierstrass representation:
\begin{equation}\label{Wcanw}
\Phi: \quad
\begin{array}{rlr}
\phi_1 &=& \ds\frac{\ii}{\sqrt{w'_1 w'_2}} \cosh \ds\frac{w_1+w_2}{2},\\[6mm]
\phi_2 &=& \ds\frac{1}{\sqrt{w'_1 w'_2}} \sinh \ds\frac{w_1+w_2}{2},\\[6mm]
\phi_3 &=& \ds\frac{1}{\sqrt{w'_1 w'_2}} \cosh \ds\frac{w_1-w_2}{2},\\[6mm]
\phi_4 &=& \ds\frac{1}{\sqrt{w'_1 w'_2}} \sinh \ds\frac{w_1-w_2}{2}.\\
\end{array}
\end{equation}
The functions $(w_1,w_2)$ in this representation satisfy the following conditions
\begin{equation}\label{Wcanw_cond}
w'_1 w'_2 \neq 0; \quad w_1 + \bar w_2 \neq (2k+1)\pi\ii;\ k\in {\mathbb Z}.
\end{equation}
Conversely, if $(w_1,w_2)$ is a pair of holomorphic functions
satisfying the conditions \eqref{Wcanw_cond}, then formulas
\eqref{Wcanw} generate a minimal space-like surface of general
type, parameterized by canonical coordinates of the first type.
\end{thm}

In terms of the functions $g_1, g_2$ the canonical Weierstrass representation of minimal
space-like surfaces of general type has the following form \cite{G-K-3}:

\begin{thm}
Any minimal space-like surface ${\M}$ of general type,
parameterized by canonical coordinates of the first type, has the
following Weierstrass representation:
\begin{equation}\label{Wcang}
\Phi: \quad
\begin{array}{rll}
\phi_1 &=& \ds\frac{\ii}{2}\; \ds\frac {g_1 g_2+1}{\sqrt{g'_1 g'_2}}\,,\\[6mm]
\phi_2 &=& \ds\frac{1}{2}\; \ds\frac {g_1 g_2-1}{\sqrt{g'_1 g'_2}}\,,\\[6mm]
\phi_3 &=& \ds\frac{1}{2}\; \ds\frac {g_1 + g_2}{\sqrt{g'_1 g'_2}}\,,\\[6mm]
\phi_4 &=& \ds\frac{1}{2}\; \ds\frac {g_1 - g_2}{\sqrt{g'_1 g'_2}}\,,\\
\end{array}
\end{equation}
where $(g_1,g_2)$ is a pair of holomorphic functions satisfying
the conditions:
\begin{equation}\label{Wcang_cond}
g'_1 g'_2 \neq 0; \quad g_1 \bar g_2 \neq -1.
\end{equation}

Conversely, if $(g_1,g_2)$ is a pair of holomorphic functions
satisfying the conditions \eqref{Wcang_cond}, then formulas
\eqref{Wcang} generate a minimal space-like surface of general
type, parameterized by canonical coordinates of the first type.
\end{thm}

\subsection{Proof of Theorem 1} Let $(K,\varkappa)$ be a solution to the system \eqref{Mink-SL-NatEq}.
Then there exists a minimal space-like surface of general type, whose Gauss and normal curvatures are exactly
the functions $K$ and $\varkappa$, respectively. This surface has a canonical Weierstrass representation
of the type \eqref{Wcang}. Applying formula  \eqref{K+ikappa_Can_g12},
we obtain
\begin{equation}
\begin{array}{lll}
        K         &=& |\alpha|\Re\alpha \\
        \varkappa &=& |\alpha|\Im\alpha
\end{array}; \quad
\alpha = \ds\frac{-4g'_1 \bar {g'_2}}{(1 + g_1 \bar g_2)^2}\;,
\end{equation}
where $(g_1,g_2)$ are a pair of holomorphic functions, satisfying the conditions $g'_1 g'_2 \neq 0$
and $g_1 \bar g_2 \neq -1$.

Conversely, let the pair of holomorphic functions $g_1, g_2$ satisfy the conditions \eqref{g_1g_2-cond}.
These two functions generate by means of \eqref{Wcang} a minimal space-like
surface of general type parameterized by canonical coordinates.
The curvatures $K$ and $\varkappa$ of this surface from one side satisfy the system \eqref{Mink-SL-NatEq},
from the other side they satisfy equality \eqref{K+ikappa_Can_g12}. Therefore the formulas
\eqref{Mink-SL-NatEq-Sol-g12} generate a solution to the system
\eqref{Mink-SL-NatEq}. \prfend

The proof of Theorem 1 is based on Theorem 2.3. Next we give another forms of Theorem 1 based on
Theorem 2.1 and Theorem 2.2.

Applying the canonical Weierstrass representation of the type \eqref{Wcanh} and formulas
\eqref{K+ikappa_Can_theta}, we get the following statement.

{\bf Theorem 1a.} {\it If the pair $(K,\varkappa)$ is a solution to the system \eqref{Mink-SL-NatEq}, then
the curvatures $K$ and $\varkappa$ are given locally by the formulas
\begin{equation}\label{Mink-SL-NatEq-Sol-theta}
\begin{array}{lll}
        K         &=& |\alpha|\Re\alpha \\
        \varkappa &=& |\alpha|\Im\alpha
\end{array}; \quad
\alpha = \ds\frac{-({\theta'_u}^2+{\theta'_v}^2)}{\cosh^2\theta}\;,
\end{equation}
where $\theta (u,v)$ is a complex harmonic function, satisfying the conditions $\cosh\theta \neq 0$
and ${\theta'_u}^2+{\theta'_v}^2 \neq 0$.

Conversely, any complex harmonic function $\theta$, satisfying the conditions $\cosh\theta \neq 0$
and ${\theta'_u}^2+{\theta'_v}^2 \neq 0$ generates a solution to the system \eqref{Mink-SL-NatEq}}
\vskip 2mm
Further, applying the canonical Weierstrass representation of the type \eqref{Wcanw}
and formulas \eqref{K+ikappa_Can_w12}, we have

{\bf Theorem 1b.} {\it If the pair $(K,\varkappa)$ is a solution to the system \eqref{Mink-SL-NatEq}, then
the curvatures $K$ and $\varkappa$ are given locally by the formulas

\begin{equation}\label{Mink-SL-NatEq-Sol-w12}
\begin{array}{lll}
        K         &=& |\alpha|\Re\alpha \\
        \varkappa &=& |\alpha|\Im\alpha
\end{array}; \quad
\alpha = \ds\frac{-w'_1 \bar {w'_2}}{\cosh^2\frac{w_1 + \bar
w_2}{2}}\;,
\end{equation}
where $(w_1,w_2)$ are holomorphic functions, satisfying the conditions $\cosh^2\frac{w_1 + \bar w_2}{2} \neq 0$
and $w'_1 {w'_2} \neq 0$.

Conversely, any such a pair of holomorphic functions generates by means of the formulas \eqref{Mink-SL-NatEq-Sol-w12}
a solution to the system \eqref{Mink-SL-NatEq}.}

\subsection {Proof of Theorem 2.} Let us denote by $\M$ and $\hat\M$ the minimal space-like surfaces of general type
given by \eqref{Wcang}. Then the condition that $(g_1,g_2)$ and $(\hat g_1,\hat g_2)$ generate one and the same
solution to \eqref{Mink-SL-NatEq} means that $\M$ and $\hat\M$ have the same curvatures $K$ and $\varkappa$.
This is equivalent with the condition that $\M$ and $\hat\M$ are obtained one from the other by a rotation
from the group $\mathbf{SO}(3,1,\RR)$ in $\RR^4_1$.

In \cite{G-K-2} we proved the following statement.
\begin{thm}\label{A-B_Cang}
Let $({\hat\M},\hat\x)$ and $({\M},\x)$ be two minimal space-like surfaces of general type, given by
the canonical Weierstrass representation of the type \eqref{Wcang}. The following conditions are equivalent:
\begin{enumerate}
    \item $({\hat\M},\hat\x)$ and $({\M},\x)$ are related by a transformation in $\RR^4_1$ of the type:\\
    $\hat\x(t)=A\x(t)+\vb$, where $A \in \mathbf{SO}(3,1,\RR)$ and $\vb \in \RR^4_1$.
    \item The functions in the Weierstrass representations of $({\hat\M},\hat\x)$ and $({\M},\x)$ are related by the
    following equalities:
\begin{equation}\label{hatg12_g12}
\hat g_1 = \ds\frac{ag_1+b}{cg_1+d}\,; \quad \hat g_2 = \ds\frac{\phantom{-}\bar d g_2 - \bar c}{-\bar b g_2 + \bar a}\;,
\end{equation}
    where $a,b,c,d\in \CC$, $ad-bc\neq 0$.
\end{enumerate}
\end{thm}

Applying theorem \ref{A-B_Cang}, we obtain the assertion. \prfend

Finally we show how to derive \eqref{Mink-SL-NatEq-Sol-1} from formula \eqref{Mink-SL-NatEq-Sol-g12}.

Let us introduce the following pair $(\xi_1,\xi_2)$ of holomorphic functions:
\[
g_1=\frac{1}{\xi_1}; \quad g_2=\xi_2\,.
\]
Then, by virtue of \eqref{Mink-SL-NatEq-Sol-g12} we get the following formula for $\alpha$:
\begin{equation}\label{alpha_xi12}
\alpha = \ds\frac{4\xi'_1 \bar {\xi'_2}}{{(\xi_1 + \bar \xi_2)}^2}\;,
\end{equation}
where $(\xi_1,\xi_2)$ is a pair of holomorphic functions satisfying the conditions $\xi_1+\bar \xi_2 \neq 0$
and $\xi'_1 {\xi'_2} \neq 0$. Further we find a complex harmonic function $\eta$,
which is related to the pair $(\xi_1,\xi_2)$ in the same way as $\theta$ is related to the pair $(w_1,w_2)$.
Then $\eta$ and $(\xi_1,\xi_2)$ satisfy the same formulas as $\theta$ and $(w_1,w_2)$.

Especially we have
\begin{equation}\label{eta_xi12}
\eta =\frac{\xi_1 + \bar \xi_2}{2}; \quad {\eta'_u}^2+{\eta'_v}^2=\xi'_1 \bar {\xi'_2}\;.
\end{equation}
Conversely, if $\eta$ is a complex harmonic function, we can obtain $(\xi_1,\xi_2)$ in the same way as
$\theta$ gives $(w_1,w_2)$. Then \eqref{eta_xi12} will be valid. Replacing $(\xi_1,\xi_2)$ by means of
$\eta$ in \eqref{alpha_xi12} we find:

{\bf Theorem 1c.} {If the pair $(K,\varkappa)$ is a solution to the system \eqref{Mink-SL-NatEq}, then $K$ and $\varkappa$
are given by:
\begin{equation}\label{Mink-SL-NatEq-Sol-eta}
\begin{array}{lll}
        K         &=& |\alpha|\Re\alpha \\
        \varkappa &=& |\alpha|\Im\alpha
\end{array}; \quad
\alpha = \ds\frac{{\eta'_u}^2+{\eta'_v}^2}{\eta^2}\;,
\end{equation}
where $\eta (u,v)$ is a complex harmonic function, satisfying the conditions $\eta \neq 0$
and ${\eta'_u}^2+{\eta'_v}^2 \neq 0$.

Conversely, any such a function $\eta$ gives by means of the formulas \eqref{Mink-SL-NatEq-Sol-eta}
a solution to the system \eqref{Mink-SL-NatEq}.}

Thus we obtained that \eqref{Mink-SL-NatEq-Sol-1} describes all solutions to the system \eqref{Mink-SL-NatEq}.

\end{document}